\date{}
\documentclass [11pt]{article}
\usepackage{amsfonts,amsmath}

\title{Shadows and intersections: stability and new proofs}
\author{
Peter Keevash \thanks{
School of Mathematical Sciences,
Queen Mary, University of London,
Mile End Road, London,
E1 4NS, UK. \newline
Email: p.keevash@qmul.ac.uk.
Research supported in part by NSF grant DMS-0555755.}
}

\newtheorem{theo}{Theorem} 

\newtheorem{lemma}[theo]{Lemma}

\newcommand{\mb}[1]{\mathbb{#1}}
\newcommand{\nib}[1]{\noindent {\bf #1}}

\newcommand{\sm}{\setminus}

\def\qed{\hfill $\Box$}

\begin{document}
\maketitle

\begin{abstract}
  We give a short new proof of a version of the Kruskal-Katona theorem
  due to Lov\'asz. Our method can be extended to a stability result,
  describing the approximate structure of configurations that are
  close to being extremal, which answers a question of Mubayi.  This
  in turn leads to another combinatorial proof of a stability theorem
  for intersecting families, which was originally obtained by Friedgut
  using spectral techniques and then sharpened by Keevash and Mubayi
  by means of a purely combinatorial result of Frankl.  We also give
  an algebraic perspective on these problems, giving yet another proof
  of intersection stability that relies on expansion of a certain
  Cayley graph of the symmetric group, and an algebraic generalisation
  of Lov\'asz's theorem that answers a question of Frankl and
  Tokushige.
\end{abstract}

\section{Introduction}

The Kruskal-Katona theorem \cite{Kat,Kr}
is a classical result in Extremal Combinatorics
that gives a tight lower bound on the size of the shadow of a
$k$-graph.\footnote{A $k$-graph $G$ consists of a
vertex set $V(G)$ and an edge set
$E(G)$, each edge being some $k$-tuple of vertices. Its shadow
$\partial G$ is the $(k-1)$-graph consisting of all $(k-1)$-tuples
that are contained in some edge of $G$. We write $|G|=|E(G)|$ for
the number of edges in $G$.}
It states that $|\partial G| \ge |\partial G_0|$, where
$G_0$ is the initial segment of length $|G|$
in the colexicographic order\footnote{If $(X,<)$ is an ordered set
we order subsets of $X$ by $A < B$ iff the largest element of
$(A \cup B) \sm (A \cap B)$ lies in $B$.} on $k$-tuples of
some ordered set. The quantative form of this statement
is a bit technical, and it is often more convenient to use the following
slightly weaker form due to Lov\'asz \cite{Lo} Ex 13.31(b): if
$|G| = \binom{x}{k} = x(x-1) \cdots (x-k+1)/k!$ for some
real number $x \ge k$ then $|\partial G| \ge \binom{x}{k-1}$.
He also showed that equality occurs if and only if $x$ is an integer
and $G = K^k_x$ is the complete $k$-graph on $x$ vertices.

This result has many consequences in Extremal Combinatorics
(see \cite{FT}). Also, its isoperimetric nature leads to broader
applications, such as the proof of the existence of threshold functions
for monotone properties by Bollob\'as and Thomason \cite{BT}.
It can also be interpreted as giving an upper bound on the number
of copies of $K^r_{r+1}$ in an $r$-graph $G$ in terms of $|G|$
(setting $r=k-1$). The general
question of estimating the number of copies of one hypergraph in
another was studied in \cite{Al} and \cite{FK}. The latter paper
gives two general bounds, one using Shearer's entropy inequality and
another using the Bonami-Beckner hypercontractive estimate. These
bounds give the correct order of magnitude in many interesting
cases, but fall short of giving the correct constant of proportionality
for complete $r$-graphs.

There are many known proofs of the Kruskal-Katona theorem
(see \cite{Da1,Hi,Fr1,Lo}) relying on compression techniques and/or
induction arguments. We start by giving a new proof (not using
either of these methods) of an upper bound on
$K^r_{r+1}(G)$, the number
of copies of $K^r_{r+1}$ in an $r$-graph $G$,
in terms of $|G|$. This can be easily translated into Lov\'asz's
result by noting that if $G$ is a $k$-graph then
$G \subset K^k_{k-1}(\partial G)$. Our proof
has the advantages that it is very simple, and the idea
can be used to obtain certain structural information not available
with other arguments.

\begin{theo} \label{kk} {\bf (Lov\'asz \cite{Lo})}
Suppose $r \ge 1$ and
$G$ is a $r$-graph with $\binom{x}{r}$ edges, for some real number
$x \ge r$. Then $K^r_{r+1}(G) \le \binom{x}{r+1}$, with equality if and only
if $x$ is an integer and $G = K^r_x$.
\end{theo}

Building on the idea in our proof of Theorem \ref{kk}, we can describe the
approximate structure of an $r$-graph $G$ that is close to being
extremal. We show that shadows have `stability',
a phenomenon which was originally discovered by Erd\H{o}s and Simonovits
in the 60's in the context of graphs with excluded subgraphs,
but has only been systematically explored relatively recently, as
researchers have realised the importance and applications of such
results in hypergraph Tur\'an theory, enumeration of discrete
and extremal set theory (see \cite{KM} as an example and for further
references). Answering a question of Mubayi (personal communication)
we prove the following stability version of the Kruskal-Katona
theorem.

\begin{theo} \label{kkstab-weak}
For any $\epsilon>0$ and $r \ge 1$ there is $\delta>0$
so that if $G$ is an $r$-graph with $\binom{x}{r}$ edges
and $K^r_{r+1}(G) > (1-\delta) \binom{x}{r+1}$ then
there is a set $S$ of $\lceil x \rceil$ vertices
so that all but at most $\epsilon \binom{x}{r}$ edges
of $G$ are contained in $S$.
\end{theo}

In fact, we can obtain further structural information and
quantify the dependance of $\delta$ on $r$ and $\epsilon$
to sufficient precision to deduce a stability theorem for
intersecting $r$-graphs. An $r$-graph is said to be intersecting
if every two of its edges have at least one common vertex.
A classical theorem of Erd\H{o}s, Ko and Rado \cite{EKR} states
that an intersecting $r$-graph $G$ on $n \ge 2r$ vertices\footnote{
The case $n<2r$ is trivial, as then $K^r_n$ is intersecting.} has
at most $\binom{n-1}{r-1}$ edges, and for $n>2r$ equality holds only when
there is some vertex $x$ that belongs to every edge of $G$.
Using spectral techniques,
Friedgut \cite{Fri} obtained a stability version, namely that, given
$\zeta>0$ there is $c>0$ so that if $\zeta n < r < (1/2-\zeta)n$
and $G$ is an intersecting $r$-graph on $n$ vertices with
$|G| > (1-\delta) \binom{n-1}{r-1}$, for some $\delta>0$,
then there is some vertex $x$ that belongs to all but at most
$c\delta \binom{n}{r}$ edges of $G$. The assumption that
$r > \zeta n$ was eliminated by Dinur and Friedgut \cite{DF}.
With $r = n/2 - t$ and $0 < t=o(n)$ one needs a lower bound of
$|G| > (1 - O(t/n)) \binom{n-1}{r-1})$ for a stability result
to hold, and such a result was obtained by Keevash and Mubayi \cite{KM}
using a purely combinatorial result of Frankl \cite{Fr2}.
Frankl's argument relies heavily on compression techniques,
but our methods give a direct proof of the following theorem,
which although weaker than that in \cite{KM} gives structural
information for all $r<n/2$.

\begin{theo} \label{intstab1}
Suppose $1 \le r < n/2$, $\delta < 10^{-3}n^{-4}$ and
$G$ is an intersecting $r$-graph on $n$ vertices
with $|G| > (1-\delta)\binom{n-1}{r-1}$. Then there is
some vertex $v$ so that all but at most
$25n\delta^{1/2} \binom{n-1}{r-1}$ edges of $G$ contain $v$.
\end{theo}

Next we take an algebraic perspective on the problem and give
yet another proof of stability, this time using expansion of
a certain Cayley graph of the symmetric group $S_{n-1}$.\footnote{
There is no similarity with the methods in \cite{Fri} and \cite{DF}
which use Fourier analysis on $\mb{Z}_2^n$.}
Here we need to assume a stronger lower bound on $|G|$, but the method
seems interesting in its own right, and has potential applications
to other problems.

\begin{theo} \label{intstab2}
Suppose $1 \le r < n/2$, $\delta <  \frac{1}{2rn^4}$ and
$G$ is an intersecting $r$-graph on $n$ vertices
with $|G| > (1-\delta)\binom{n-1}{r-1}$. Then there is
some vertex $v$ so that all but at most
$\delta r \binom{n-1}{r-1}$ edges of $G$ contain $v$.
\end{theo}

Given an $r$-graph $G$ there are some naturally associated
algebraic objects called (higher) inclusion matrices.
For $s \le r$ we define $M^r_s(G)$ as a $\{0,1\}$ matrix with
rows indexed by edges of $G$ and columns indexed by subsets of $V(G)$
of size $s$: the entry corresponding to an edge $e$ and a set $S$
is $1$ if $S \subset e$ and $0$ otherwise. Frankl and Tokushige \cite{FT}
posed the problem of determining the minimum rank
$\mbox{rk } M^r_s(G)$ of $M^r_s(G)$
in terms of $|G|$. We obtain the following result.

\begin{theo} \label{kk-alg}
For every $r \ge s \ge 0$ there is a number $n_{r,s}$ so that
if $G$ is an $r$-graph with $|G| = \binom{x}{r} \ge n_{r,s}$ then
$\mbox{rk } M^r_s(G) \ge \binom{x}{s}$.
If $r>s>0$ then equality holds only
if $x$ is an integer and $G=K^r_x$.
\end{theo}

Note that this generalises the result of Lov\'asz, and also its iterated
version, i.e. that if $G$ is an $r$-graph, $|G| = \binom{x}{r}$
and $s \le r$ then $|\partial^r_s G| \ge \binom{x}{s}$, where
the $s$-shadow $\partial^r_s G$ consists of all $s$-sets
that are contained in some edge of $G$. This is immediate
from Theorem \ref{kk-alg} (for large $x$) since the rank
of $M^r_s(G)$ is at most the number of non-zero columns, which
is the size of the $s$-shadow. Keevash and Sudakov \cite{KS}
obtained a non-uniform version of this inequality, and our proof
uses elements of that approach, but requires a number of new
ideas. We highlight one lemma that we think is of independent
interest, as it expresses a certain rigidity property of
the complete inclusion matrices $M^r_s(K^r_n)$.

\begin{lemma} \label{full-rk}
Suppose $0 \le s \le r < n/2$ and $G = K^r_n \sm F$ is an $r$-graph
on $[n]$ with $|F| < \binom{r}{s}^{-1} \binom{n}{r-s}$. Then
$\mbox{rk } M^r_s(G) = \binom{n}{s}$.
\end{lemma}

The rest of this paper is organised as follows. The next section
gives a very short proof of Theorem \ref{kk}. In section 3 we
collect some facts about binomial coefficients and other
inequalities that will be
subsequently useful. In section 4 we extend the ideas from our
proof of Theorem \ref{kk} to prove a generalised form of
Theorem \ref{kkstab-weak}. This is then combined with an idea
of Daykin in the following section to obtain our first proof
of stability for intersecting families. Section 6 contains
our second proof, based on expansion in the symmetric group.
In section 7 we prove our bound on the rank of inclusion matrices,
Theorem \ref{kk-alg}, and the final section contains some concluding
remarks.

\nib{Notation.} We write $[n]=\{1,\cdots,n\}$.
Suppose $G$ is an $r$-graph. Let $K^r_{r+1}(G)$ be the number of copies
of $K^r_{r+1}$ in $G$. For a vertex $v \in V(G)$ let
$K^r_{r+1}(v)$ be the number of $K^r_{r+1}$'s that contain $v$.
The link $(r-1)$-graph is
$L(v) = \{A \subset V(G): |A|=r-1, A \cup \{v\} \in E(G)\}$.
The degree $d(v)=|L(v)|$ is the number of edges containing $v$.

\section{Proof of Theorem \ref{kk} }

We argue by induction on $r$. The base case $r=1$ is trivial.
We can assume that the degree $d(v)$ is non-zero for every vertex $v$.
Note that $S \cup \{v\}$ spans a $K^r_{r+1}$ in $G$ if and only if
$S$ is an edge of $G$ and spans a $K^{r-1}_r$ in the link
$L(v)$. The first
condition gives the estimate $K^r_{r+1}(v) \le |G|-d(v)$
and the second $K^r_{r+1}(v) \le K^{r-1}_r(L(v))$.
We claim that $K^r_{r+1}(v) \le (x/r-1)d(v)$ for every $v$, and equality
is only possible when $d(v) = \binom{x-1}{r-1}$.
To see this, suppose first that $d(v) \ge \binom{x-1}{r-1}$. Then
by the first condition it suffices to observe that
$\binom{x}{r} - d(v) \le (x/r-1)d(v)$. On the other hand, if
$d(v) \le \binom{x-1}{r-1}$ then
define the real number $x_v \ge r$ by $d(v) = \binom{x_v-1}{r-1}$.
Then by induction hypothesis
$K^{r-1}_r(L(v)) \le \binom{x_v-1}{r} = (x_v/r-1)d(v) \le (x/r-1)d(v)$.
The equality conditions are clear, so the claim holds in either case.
Now
$$(r+1)K^r_{r+1}(G) = \sum_v K^r_{r+1}(v) \le (x/r-1)\sum_v d(v)
= (x/r-1)r|G| = (x-r)\binom{x}{r} = (r+1)\binom{x}{r+1}.$$
Therefore $K^r_{r+1}(G) \le \binom{x}{r+1}$, as required. Equality holds
only when all vertices have degree $\binom{x-1}{r-1}$. Then
if $G$ has $n$ vertices we have
$n\binom{x-1}{r-1} = \sum_v d(v) = r \binom{x}{r}$, so $n=x$ and
$G=K^r_x$. \qed

\section{Technical estimates}

We pause to collect some technical estimates
that will be helpful in the following sections.
The first two concern binomial coefficients, and we will
prove them in an appendix to the paper. The others
are straightforward, so we omit the proofs. We consider
the binomial coefficient $\binom{x}{r}$ to be the polynomial
$x(x-1)\cdots(x-r+1)/r!$ defined for every real number $x$.
It is positive and increasing for $x>r-1$.

\begin{lemma} \label{bin-diff}
If $x>y \ge r-1$ then $(x-y) \binom{y-1}{r-1} < \binom{x}{r} - \binom{y}{r}
 < (x-y) \binom{x}{r-1}$.
\end{lemma}

\begin{lemma} \label{bin-shadow}
Suppose $r \ge 2$, $1 \le s \le r-1$,
$\binom{u}{r} = \binom{v}{r} + \binom{w}{r-1}$
with $1 \le \binom{w}{r-1} < \binom{u-1}{r-1} - \frac{1}{2r!} u^{r-s-1}$,
and $u > u_0(r,s)$ is sufficiently large.
Then $\binom{u}{s} < \binom{v}{s} + \binom{w}{s-1}
- (3r)^{-r}u^{-1}$.
\end{lemma}

There follow some assorted easy facts. Throughout
$n$ is a natural number, other parameters are real.

\begin{equation} \label{f1}
\mbox{If } 0<\theta <1 \mbox{ and } \theta x > n
\mbox{ then } \binom{\theta x}{n} < \theta^n \binom{x}{n}.
\end{equation}

\begin{equation} \label{f2}
\mbox{If } a>b>0 \mbox{ then }
\left( \frac{a}{b} \right)^n < \binom{a}{n}/\binom{b}{n}
= \prod_{i=0}^{n-1} \frac{a-i}{b-i} < \left( \frac{a-n+1}{b-n+1} \right)^n.
\end{equation}

\begin{equation} \label{f3}
\mbox{If } 0 < \theta < 2/3n \mbox{ and }
\mbox{ then } (1+\theta)^n < 1 + 2n\theta.
\end{equation}

\begin{equation} \label{f4}
\mbox{If } 0 < \theta < 1/2n \mbox{ then } (1+\theta)^n < 2.
\end{equation}

\begin{equation} \label{f5}
\mbox{If } 0 < \theta < 1/2 \mbox{ then } (1-\theta)^{1/n} > 1-2\theta/n.
\end{equation}

\section{Stability for shadows}

Building on the idea in our proof of Theorem \ref{kk}, we can describe the
approximate structure of an $r$-graph $G$ that is close to being
extremal. Answering a question of Mubayi (personal communication)
we obtain a quantative stability version of the Kruskal-Katona
theorem: statement (1) in the following theorem.

\begin{theo} \label{kkstab}
Suppose $0 < \epsilon < 1/2$, $r \ge 2$, $x \ge (1+\epsilon)(r+1)$,
$\delta < (\epsilon/6r)^2$ and $G$ is an $r$-graph with
$\binom{x}{r}$ edges and
$K^r_{r+1}(G) > (1-\delta) \binom{x}{r+1}$. Then:

\noindent (1) There is a set $S$ of $\lceil x \rceil$ vertices
so that all but at most $10r(\epsilon^{-1}+1)\delta^{1/2} \binom{x}{r}$
edges of $G$ are contained in $S$.

\noindent (2) There are at most $\delta^{1/2}x$ vertices
$v$ with $d(v) >  (1+\delta^{1/2}) \binom{x-1}{r-1}$.

\noindent (3) The vertices of degree less than
$\binom{(1-\delta^{1/2})(x-1)}{r-1}$ are incident to at most
$\delta^{1/2} x \binom{x-1}{r-1}$ edges.

\noindent (4) There is a set $C$ of size $|C| < (1+3r\delta^{1/2})x$
that contains at least $\binom{(1-4\delta^{1/2})x}{r+1}$ copies
of $K^r_{r+1}$.
\end{theo}

\nib{Proof.} Note that our assumption $K^r_{r+1}(G)>0$ implies
that $x \ge r+1$. For each vertex $v$ we recall the
bounds $K^r_{r+1}(v) \le \binom{x}{r}-d(v)$,
and $K^r_{r+1}(v) \le (x/r-1)d(v)$ proved above, and
the bound $K^r_{r+1}(v) \le \binom{x-1}{r}$, which follows
by combining the first two bounds. Also, if $d(v) = \binom{x_v-1}{r-1}$
we recall that $K^r_{r+1}(v) \le (x_v/r-1)d(v)$.
Let
$$A = \left\{v: d(v) > \binom{x-1}{r-1} \right\}, \qquad
A_0 = \left\{v: d(v) > (1+\delta^{1/2}) \binom{x-1}{r-1} \right\},$$
$$B = V(G) \sm A = \left\{v: d(v) \le \binom{x-1}{r-1} \right\},
\qquad \mbox{ and }
B_0 = \left\{v: d(v) < \binom{x-1-y}{r-1} \right\},$$
where $y = \delta^{1/2}(x-r)$.
For a set $S$ write $d_S = \sum_{v \in S} d(v)$.
We have
\begin{eqnarray*}
(1-\delta)(r+1) \binom{x}{r+1} &<& (r+1) K^r_{r+1}(G) = \sum_v K^r_{r+1}(v)
= \sum_{v \in A} K^r_{r+1}(v) + \sum_{v \in B} K^r_{r+1}(v) \\
& \le & |A|\binom{x}{r}-d_A + (x/r-1)d_B\\
&=& |A|\binom{x}{r} + (x/r-1)(d_A+d_B) - (x/r)d_A
\end{eqnarray*}
and
$(x/r-1)(d_A+d_B) = (x-r)\binom{x}{r}=(r+1)\binom{x}{r+1}$,
so $|A| > d_A \binom{x-1}{r-1}^{-1} - \delta(x-r)$.
Now $$|A_0| \delta^{1/2} \binom{x-1}{r-1} + |A| \binom{x-1}{r-1}
< d_A <  \binom{x-1}{r-1} (|A|+\delta(x-r)),$$
so $|A_0| < \delta^{1/2}(x-r)$. This implies (2). We also deduce
\begin{equation} \label{a0}
\sum_{v \in A_0} K^r_{r+1}(v) \le  |A_0| \binom{x-1}{r} <
\delta^{1/2}(x-r)\binom{x-1}{r}
< \delta^{1/2}(r+1)\binom{x}{r+1}.
\end{equation}
Next we have
\begin{eqnarray*}
(1-\delta)(r+1) \binom{x}{r+1} &<& (r+1) K^r_{r+1}(G) = \sum_v K^r_{r+1}(v)
\le (x/r-1)d_A + \sum_{v \in B} (x_v/r-1)d(v)\\
&<& (x/r-1)(d_A+d_B) - (y/r) d_{B_0} =
(r+1)\binom{x}{r+1} - (y/r) d_{B_0},
\end{eqnarray*}
so $d_{B_0} < \delta^{1/2} r \binom{x}{r}$. This implies (3). We
also deduce
\begin{equation} \label{b0}
\sum_{v \in B_0} K^r_{r+1}(v) \le (x/r-1)d_{B_0} <
\delta^{1/2}(r+1)\binom{x}{r+1}.
\end{equation}

Define an $(r+1)$-graph $H$ on the same vertex set of $G$ where
an $(r+1)$-tuple is an edge exactly when it spans a $K^r_{r+1}$ in $G$.
Let $C = V(G) \sm (A_0 \cup B_0)$,
$H_0 \subset H$ consist of all $(r+1)$-tuples of $H$ that are
contained in $C$ and $H_1 = H \sm H_0$.
Using equations (\ref{a0}) and (\ref{b0}) we have
\begin{eqnarray*}
 (1-3\delta^{1/2})(r+1)\binom{x}{r+1} & < & \sum_{v \in C} K^r_{r+1}(v)
= \sum_{v \in C} d_H(v) < (r+1)|H_0| + r|H_1| \\
& = & (r+1)|E(H)|-|H_1| < (r+1)\binom{x}{r+1} - |H_1|,
\end{eqnarray*}
where in the last step we use Theorem \ref{kk}. Therefore
$|H_1| < 3\delta^{1/2}(r+1)\binom{x}{r+1}$ and
\begin{eqnarray}
|H_0| & > & (r+1)^{-1}\left( (1-3\delta^{1/2})(r+1)\binom{x}{r+1}
- r|H_1| \right) \nonumber \\
& > & (1-3(r+1)\delta^{1/2})\binom{x}{r+1}
> \binom{(1-4\delta^{1/2})x}{r+1}, \label{h0}
\end{eqnarray}
where in the last step we use fact (\ref{f1}) to obtain
the inequality
$$\binom{x}{r+1}^{-1}\binom{(1-4\delta^{1/2})x}{r+1}
< (1-4\delta^{1/2})^{r+1} < 1 - (r+1) \cdot 4\delta^{1/2}
+ \binom{r+1}{2} (4\delta^{1/2})^2 < 1 - 3(r+1)\delta^{1/2},$$
which is valid since $\delta < (\epsilon/6r)^2 < 1/64r^2$.
Now we can apply Theorem \ref{kk} to $H_0$ to deduce that
at least $\binom{(1-4\delta^{1/2})x}{r}$ edges of $G$ are
contained in $C$.

Next, since
$r \binom{x}{r} > \sum_{v \in C} d(v) \ge |C| \binom{x-1-y}{r-1}$
by fact (\ref{f2}) we have
$$|C|/x \le \binom{x-1}{r-1}\binom{x-y-1}{r-1}^{-1}
< \left( \frac{x-r+1}{x-y-r+1} \right)^{r-1}
< 1 + \frac{2(r-1)\delta^{1/2}}{1-\delta^{1/2}},$$
where we use fact (\ref{f3}).
Thus we can write
$|C|=x+t$ with $t/x < 2(r-1)\delta^{1/2}(1-\delta^{1/2})^{-1}
< 3r\delta^{1/2}$, which proves (4).

Choose any set $S \subset V(G)$ with
$|S| = \lceil x \rceil$ so that either $S \subset C$ if
$|C| \ge \lceil x \rceil$ or $S \supset C$ if
$|C| \le \lceil x \rceil$. We can estimate the number of edges
of $G$ that are not contained in $S$ as follows.
Either such an edge is not contained in $C$, of which there
are at most
$$|G|-|H_0| < \binom{x}{r} -\binom{(1-4\delta^{1/2})x}{r}
< 4\delta^{1/2}x \binom{x}{r-1}
= \frac{4r\delta^{1/2}}{1-(r+1)/x} \binom{x}{r} <
4(\epsilon^{-1}+1)r\delta^{1/2}\binom{x}{r},$$
or it is contained in $C$ but not in $S$, of which there are at most
\begin{eqnarray*}
\binom{|C|}{r} - \binom{|S|}{r} & < & t \binom{x+t}{r-1}
< t \left( \frac{x+t-r+1}{x-r+1} \right)^{r-1} \binom{x}{r-1} \\
& = & \frac{rt}{x-r+1} \left( 1 + \frac{t}{x-r+1} \right)^{r-1}
 \binom{x}{r} < 6r(\epsilon^{-1}+1) \delta^{1/2}\binom{x}{r}.
\end{eqnarray*}
(In both estimates we use Lemma \ref{bin-diff}.
In the last step we use the estimate
$\frac{t}{x-r+1} < \frac{3r\delta^{1/2}}{\epsilon(r+1)}
< \frac{3r \cdot \epsilon/6r}{\epsilon(r+1)} < 1/2(r-1)$
and so by fact (\ref{f4}) we have
$\left( 1 + \frac{t}{x-r+1} \right)^{r-1} < 2$.)
In total we have
at most $10r(\epsilon^{-1}+1)\delta^{1/2} \binom{x}{r}$ edges of $G$
not contained in $S$, as required. \qed

\medskip

\nib{Remark.} We have tried to give good estimates in this proof so
that we obtain stability results for a large range of $r$ and $x$, but
some price has been paid for obtaining a universal bound,
and improvements can be made for particular values of the parameters.
The bounds get worse for smaller $x$ to the point where we
lose an exponential factor in $r$ if $x=r+c$ and $c \ll r$.
The proof breaks down as $x$ approaches $r+1$, but in this range
the weak Kruskal-Katona bound compares poorly to the full theorem,
and in any case it is not too hard to analyse the situation
by ad hoc methods.

\section{Stability for intersecting familes, I}

Next we show how to derive a stability result for intersecting families.
The proof involves combining the methods above with
an idea of Daykin \cite{Da2}.
First we remark that if $r \le l \le m$ and $G$ is a $r$-graph
with  $K^r_l(G) = \binom{x}{l}$ then
$K^r_m(G) \le \binom{x}{m}$. This follows by repeatedly
applying Theorem \ref{kk} and noting that a set $M$ of size $m$ spans
a $K^r_m$ in $G$ exactly when it spans a $K^{m-1}_m$ in the
$(m-1)$-graph of all copies of $K^r_{m-1}$ in $G$.

Now we prove Theorem \ref{intstab1}, which is as follows:
Suppose $1 \le r < n/2$, $\delta < 10^{-3} n^{-4}$ and
$G$ is an intersecting $r$-graph on $n$ vertices
with $|G| > (1-\delta)\binom{n-1}{r-1}$. Then there is
some vertex $v$ so that all but at most
$25n\delta^{1/2}\binom{n-1}{r-1}$ edges of $G$ contain $v$.

\nib{Proof.}\
Consider the complementary $r$-graph
$H = \{A \subset V(G): |A|=r, A \notin E(G)\}$ and
the $(n-r)$-graph of complements
$J = \{A \subset V(G): V(G) \sm A \in E(G) \}$.
Write $|H| = \binom{n-\theta}{r}$.
By the Erd\H{o}s-Ko-Rado theorem
we have $|G| \le \binom{n-1}{r-1}$, so $|H| \ge \binom{n-1}{r}$,
i.e. $0 \le \theta \le 1$. Write $\phi=1-\theta$.
By Lemma \ref{bin-diff} we have
$\binom{n-\theta}{r} - \binom{n-1}{r} > \phi \binom{n-2}{r-1}$, so
$$\binom{n-1}{r} + \phi \binom{n-2}{r-1} < |H|
= \binom{n}{r}-|G| < \binom{n-1}{r} + \delta \binom{n-1}{r-1},$$
and $\phi < \delta \frac{n-1}{n-r} < 2\delta$
(since $r < n/2$).
The condition that $G$ is intersecting may be rephrased as
saying that every edge of $J$ spans a $K^r_{n-r}$ in $H$.
Therefore
\begin{eqnarray*}
K^r_{n-r}(H) & \ge & |J| =  \binom{n}{r}-\binom{n-\theta}{r}
= \left( \binom{n}{r} - \binom{n-1}{r} \right)
- \left( \binom{n-\theta}{r}-\binom{n-1}{r} \right) \\
& > & \binom{n-1}{r-1} - \phi \binom{n}{r-1}
 = \left( 1 - \frac{n}{n-r+1}\phi \right) \binom{n-1}{n-r} \\
& > & \left( 1 - \frac{n}{n-r+1}\phi \right)
   \left( \binom{n-\theta}{n-r} - \phi \binom{n-\theta}{n-r-1} \right)\\
& = & \left( 1 - \frac{n}{n-r+1}\phi \right)
     \left( 1 - \frac{n-r}{r+\phi} \phi \right) \binom{n-\theta}{n-r}
  > \left( 1 - \frac{n}{n-r+1}\phi
     - \frac{n-r}{r+\phi} \phi \right) \binom{n-\theta}{n-r}\\
& > & (1-4(n/r-1)\phi) \binom{n-\theta}{n-r}.
\end{eqnarray*}
Write $c = 1-4(n/r-1)\phi$. Then there must be some
$m$ with $r \le m < n-r$ for which
$$K^r_m(H) \le c^{\frac{m-r}{n-2r}} \binom{n-\theta}{m}
\quad \mbox{ and } \quad
K^r_{m+1}(H) \ge c^{\frac{m+1-r}{n-2r}} \binom{n-\theta}{m+1}.$$
Write $K^r_m(H) = \binom{n-\psi}{m}$, where $\psi \ge \theta$ by the
remark before the proof. Also, by the same remark we have
$$\binom{n-\psi}{n-r} \ge K^r_{n-r}(H) > c \binom{n-\theta}{n-r}
> \binom{n-2}{n-r},$$
as
$c \binom{n-\theta}{n-r} - \binom{n-2}{n-r}
> c \binom{n-1}{n-r} - \binom{n-2}{n-r}
= \left( 1 - 2\phi \frac{n}{r+1} - \frac{r-1}{n-1} \right) > 0$,
since $\phi<2\delta < 2^{-9}n^{-3}$. This gives $\psi < 2$.
Now we have
\begin{eqnarray*}
K^r_{m+1}(H) & \ge & c^{\frac{m+1-r}{n-2r}} \binom{n-\theta}{m+1}
= c^{\frac{m-r}{n-2r}} \binom{n-\theta}{m} \cdot
 c^{\frac{1}{n-2r}} \frac{n-\theta-m}{m+1}
\ge K^r_m(H) c^{\frac{1}{n-2r}} \frac{n-\theta-m}{m+1} \\
& = & c^{\frac{1}{n-2r}} \frac{n-\theta-m}{n-\psi-m} \binom{n-\psi}{m+1}
 \ge \left( 1 - \frac{8(n/r-1)\phi}{n-2r} \right) \binom{n-\psi}{m+1}
 > (1 - 25\delta) \binom{n-\psi}{m+1},
\end{eqnarray*}
where we apply fact (\ref{f5}) in the penultimate inequality and then
estimate
$\frac{8(n/r-1)\phi}{n-2r} = 8\phi(r^{-1}+(n-2r)^{-1}) < 10\phi
< 25 \delta$. \footnote{If $1 \le a \le r \le b < n/2$ then
$r^{-1}+(n-2r)^{-1}$ is maximised at $r=a$ or $r=b$, so we can
improve our bounds with more information about $r$.}
By part (4) of Theorem \ref{kkstab}, we can find a set
$C$ with $|C| = n-\psi + t$ and $t < 3m (25\delta)^{1/2}(n-\psi)
 = 15 \delta^{1/2} m(n-\psi)$
so that $H$ has at least
$\binom{(1-20\delta^{1/2})(n-\psi)}{m+1}$ copies of $K^r_{m+1}$
contained in $C$. Note that
$|C| < n - \psi + 15mn \delta^{1/2} < n$,
since $15mn \delta^{1/2} < 1/2 < 1 - 2\delta
< 1 - \phi = \theta \le \psi$. By arbitrarily adding vertices if
necessary we may assume that $|C|=n-1$.

By Theorem \ref{kk} there are at least
$\binom{(1-20\delta^{1/2})(n-\psi)}{r}$ edges of $H$ contained in $C$.
Write $\{v\} = V(G) \sm C$. Then
the number of edges of $H$ containing $v$ is at most
$Q = \binom{n-\theta}{r} - \binom{(1-20\delta^{1/2})(n-\psi)}{r}$.
Since $G$ is the complement of $H$, the number of edges of $G$
containing $v$ is at least $\binom{n-1}{r-1} - Q$,
and so by Erd\H{o}s-Ko-Rado the number of edges of $G$ not
containing $v$ is at most
\begin{eqnarray*}
|G| - \left( \binom{n-1}{r-1} - Q \right)
& \le & Q < (\phi + 20\delta^{1/2}n) \binom{n-\theta}{r-1} \\
& < & (\phi + 20\delta^{1/2}n) \left(
\binom{n-1}{r-1} + \phi \binom{n-\theta}{r-2} \right) \\
& < & (\phi + 20\delta^{1/2}n)\left( 1 + \phi \frac{n(r-1)}{(n-r+2)(n-r+1)}
\right) \binom{n-1}{r-1} \\
& = & \epsilon \binom{n-1}{r-1},
\end{eqnarray*}
where $\epsilon < (2\delta + 20\delta^{1/2}n)(1+4\delta)
< 25\delta^{1/2} n$. This completes the proof. \qed

\section{Stability for intersecting families, II}

Now we give another argument using expansion properties
of the symmetric group. We need to assume that $G$ is closer
to the maximum, but then the bound on bad edges improves.
Also, we think that the method is interesting in itself,
as it may apply to a much wider class of problems.

Our approach is based on
Katona's permutation method. We write a permutation $\sigma \in
S_n$ as a sequence $(\sigma(1), \cdots, \sigma(n))$. Say that $\sigma$
and $\tau$ are cyclically equivalent if there is some $i \in [n]$ such
that $\tau(x)=\sigma(x+i)$ for all $x \in [n]$. (Addition is mod $n$,
i.e. $x+i$ means either $x+i$ or $x+i-n$, whichever lies in $[n]$.)
Let $C_n$ be the set of equivalence classes of this relation, which
are called cyclic orders. We will abuse notation and identify a given
cyclic order with the permutation $\sigma$ that represents this class
and has $\sigma(n)=n$. Then restricting $\sigma$ to $[n-1]$
establishes a bijection between $C_n$ and $S_{n-1}$.

We consider the Cayley graph $C$ on $S_{n-1}$ generated by the set of
adjacent tranpositions $T = \{(12),(23),\cdots,(n-2\ n-1)\}$, i.e. the
vertex set of $C$ is $S_{n-1}$ and permutations $\sigma$ and $\tau$
are adjacent in $C$ if $\tau = \sigma \circ t$ for some $t \in T$.
Note that we use the multiplication convention `first $t$ then
$\sigma$', so that transpositions act by interchanging adjacent
positions (rather than values) in the sequence representing a
permutation, i.e. $(\tau(1),\cdots,\tau(n))$ is obtained
from $(\sigma(1),\cdots,\sigma(n))$ by interchanging two
consecutive elements. $C$ is a regular graph with degree $d=n-2$.
The adjacency matrix of $C$ has eigenvalues
$d = \lambda_1 \ge \lambda_2 \ge \cdots \ge \lambda_{(n-1)!}$.  A
theorem of Bacher \cite{Ba} states that the second eigenvalue
satisfies $d - \lambda_2 = 2 - 2 \cos (\pi/(n-1))$.  We will just use
the estimate $d - \lambda_2 > 2/n^2$ for $n \ge 3$, which can easily
be derived from this formula and the inequality $\cos x < 1 - x^2/4$
for $0<x<2$.

It follows that $C$ is a $\alpha$-expander, with $\alpha =
(d-\lambda_2)/2d > 1/n^3$, i.e.  for any set $W \subset V(G)$ with
$|W| \le (n-1)!/2$ we have $|N(W)| \ge |W|/n^3$, where $N(W)$ is the set of
vertices in $V(G) \sm W$ that are adjacent to some vertex of $W$.
(This value of $\alpha$ is given by Corollary 9.2 in Alon-Spencer
\cite{AS}; it is not optimal, but suffices for our purpose.)

We need the following well-known lemma, which is the basis
for Katona's proof of the Erd\H{o}s-Ko-Rado theorem. Given
a cyclic order $\sigma$, the intervals of length $r$ are the sets
$I_{\sigma,r}(x) = \{\sigma(x),\sigma(x+1),\cdots,\sigma(x+r-1)\}$
for $x \in [n]$ (addition mod $n$).

\begin{lemma} \label{cyclic}
Suppose $\sigma$ is a cyclic order of $[n]$ and
$F$ is an intersecting family of intervals of length $r<n/2$ in $\sigma$.
Then $|F| \le r$, and equality holds exactly when
there is a single point $x$ that belongs to all of the intervals.
\end{lemma}

For the convenience of the reader we include the brief proof.

\nib{Proof.}
Suppose $F$ contains the interval $I_{\sigma,r}(x)$.
Let $j \ge 0$ be maximal so that $y=x+j$ mod $n$ is in $I_{\sigma,r}(x)$
and $I_{\sigma,r}(y) \in F$. We claim
that any interval $I_{\sigma,r}(z)$ in $F$ contains $y$. To see
this, note that since $I_{\sigma,r}(z)$ intersects
$I_{\sigma,r}(x)$ we either have $z \in I_{\sigma,r}(x)$ or
$z+r-1 \in I_{\sigma,r}(x)$. In the former case we have
$z = x+j'$ mod $n$ with $0 \le j' \le j$ by definition of $y$,
so $y \in I_{\sigma,r}(z)$.
In the latter case we must have $z+r-1 = x+j'$ mod $n$
with $j \le j' \le r-1$, or
otherwise $I_{\sigma,r}(z)$ would be disjoint from
$I_{\sigma,r}(y)$, so again $y \in I_{\sigma,r}(z)$. \qed

\medskip

Now we prove Theorem \ref{intstab2}, which is as follows:
Suppose $1 \le r < n/2$, $\delta <  \frac{1}{2rn^4}$ and
$G$ is an intersecting $r$-graph on $n$ vertices
with $|G| \ge (1-\delta)\binom{n-1}{r-1}$. Then there is
some vertex $v$ so that all but at most
$\delta r \binom{n-1}{r-1}$ edges of $G$ contain $v$.

\nib{Proof.}
For each cyclic order $\sigma \in C_n$ let $G(\sigma)$ consist of
those sets of $G$ that are intervals in $\sigma$. We say
$\sigma$ is complete if $|G(\sigma)|=r$, otherwise incomplete.
The lemma tells us
that if $\sigma$ is complete then there is some point $v$
belonging to all intervals of $G(\sigma)$. To specify this point
we say that $\sigma$ is $v$-complete. Let $X$ be the set
of incomplete $\sigma$. Then
$$r!(n-r)!|G| = \sum_{\sigma \in C_n} |G(\sigma)|
\le \sum_{\sigma \in C_n \sm X} r + \sum_{\sigma \in X} (r-1)
= r(n-1)! - |X|,$$
so $|X| \le r(n-1)! - (1-\delta)r!(n-r)!\binom{n-1}{r-1}
 = \delta r(n-1)!$.
It follows that the number of complete $\sigma$ is at least
$(1 - \delta r)(n-1)!$.

Now we make the following claim:
if $\sigma$ is $v$-complete, $\tau$ is complete, and
$\tau = \sigma \circ (i\ i+1)$ for some $i \in [n] \sm \{v,v-1\}$,
then $\tau$ is $v$-complete.
To prove this, we start by relabelling (if necessary)
so that $v=n$, and so $1 \le i \le n-2$.
Since $\sigma$ is $n$-complete
we have $I_{\sigma,r}(x) \in G(\sigma)$ for
$n-r+1 \le x \le n$. Also, if $n-r+1 \le x \le n$,
$x \ne i+1$ and $x+r-1 \ne i$
then $I_{\tau,r}(x) = I_{\sigma,r}(x) \in G$ (the order is different
but the sets are the same). We have three cases according
to the value of $i$. Firstly, if $i \ne n-r$ and
$i \ne r-1$ then $I_{\tau,r}(n-r+1) = I_{\sigma,r}(n-r+1)$
and $I_{\tau,r}(n) = I_{\sigma,r}(n)$ are both in $G(\tau)$,
and their only common position is $n$, so $\tau$ must be
$n$-complete. Secondly, if $i=r-1$ then $i+1 \ne n-r+1$ (since
$r<n/2$) so $I_{\sigma,r}(x) = I_{\tau,r}(x) \in G(\tau)$
for $n-r+1 \le x \le n-1$. These intervals have just
two common positions: $n-1$ and $n$. Since $\tau$ is complete
$G(\tau)$ must either contain $I_{\tau,r}(n)$ or
$I_{\tau,r}(n-r)$. The latter case is impossible, as
$I_{\tau,r}(n-r) = I_{\sigma,r}(n-r)$ (since $i+1=r<n-r$),
but this is disjoint to $I_{\sigma,r}(n) \in G(\sigma)$
and $G$ is intersecting.
Therefore $I_{\tau,r}(n) \in G(\tau)$, i.e. $\tau$ is
$n$-complete. The argument for the third case, when $i=n-r$,
is the same as that for the second case (by symmetry),
so we will omit it. This proves the claim.

Now consider the Cayley graph $C$ on $S_{n-1}$ defined above. Suppose
$W$ is a set of complete cyclic orders, which we may consider as a
subset of $V(C)$. Since $C$ is a $1/n^3$-expander, if
$n^3 \delta r \le |W|/(n-1)! \le 1/2$ we have
$|N(W)| > \delta r (n-1)!$,
and so there is a complete $\sigma$ in $N(W)$. It
follows that the restriction of $C$ to the set of complete cyclic
orders has a connected component $C'$ of size at least
$(1 - n^3 \delta r)(n-1)!$.\footnote{
Consider the components of $C$ restricted to the complete cyclic orders.
Each component must either have size at most $n^3 \delta r (n-1)!$
(`small') or more than $(n-1)!/2$ (`large'). Since components are
disjoint sets there is at most one large component. Also, the
total size of all small components is at most $n^3 \delta r (n-1)!$,
or we could take $W$ to be a union of small components with
$n^3 \delta r (n-1)! \le |W| \le 2n^3 \delta r (n-1)!$ and
find a complete $\sigma$ in $N(W)$, contradicting the definition
of components. Therefore there is a large component, and its
size is at least $(1 - n^3 \delta r)(n-1)!$.}
By the claim, there is some $v$ so that every
$\sigma$ in $C'$ is $v$-complete.
Write $G_v$ for the sets in $G$ that contain $v$. Then
$r!(n-r)!|G_v| \ge \sum_{\sigma \in C'} |G(\sigma)| \ge
(1-\delta r)(n-1)!\cdot r$,
so $|G_v| \ge (1-\delta r)\binom{n-1}{r-1}$. Now by the Erd\H{o}s-Ko-Rado theorem
there are at most $\delta r \binom{n-1}{r-1}$ sets of $G$ that do not
contain $v$, as required. \qed

\medskip

\nib{Remark.} The generators we use in this argument are poor
from an expansion point of view, and in fact Kassabov \cite{Kas}
has shown that a constant eigenvalue gap can be obtained with
just a constant number of generators (universal constants independent
of $n$). However, this does not imply an improvement to our theorem,
as we rely heavily on structural properties of the
generating set in our argument.

\section{An algebraic generalisation of Lov\'asz's Theorem}

In this section we prove an algebraic generalisation of the
Lov\'asz version of the Kruskal-Katona theorem.
Let $G$ be an $r$-graph and $s \le r$.
The (higher) inclusion matrix $M^r_s(G)$ is a $\{0,1\}$ matrix with
rows indexed by edges of $G$ and columns indexed by subsets of $V(G)$
of size $s$: the entry corresponding to an edge $e$ and a set $S$
is $1$ if $S \subset e$ and $0$ otherwise. Frankl and Tokushige
\cite{FT} posed the problem of finding the minimum rank
of $M^r_s(G)$ in terms of $|G|$.

When $|G| = \binom{n}{r}$ for an integer $n$ then one natural
construction is the complete $r$-graph $K^r_n$. Here the
rank is given by a theorem of Gottlieb (\cite{Go}, see also
\cite{BF}):

\begin{theo} {\bf (Gottlieb \cite{Go})}
$\mbox{rk } M^r_s(K^r_n) = \min \left\{ \binom{n}{r}, \binom{n}{s}
\right\}$.
\end{theo}

Before describing what might be expected in general we describe
some recursive properties of inclusion matrices.
We define two operations associated with a vertex $x$ of $G$
giving hypergraphs on $V(G) \sm \{x\}$.
Deletion gives the $r$-graph $G \sm x = \{A: A \in G, x \notin A\}$.
Contraction gives the $(r-1)$-graph $G/x = \{A \sm \{x\}: x \in A \in G\}$.

\begin{lemma} \label{rk-induct}
Suppose $G$ is an $r$-graph, $x$ is a vertex of $G$ and $1 \le s \le r-1$.
Then
$$\mbox{rk } M^r_s(G) \ge \max\{
\mbox{rk } M^r_s(G \sm x) + \mbox{rk } M^{r-1}_{s-1}(G/x),
\mbox{rk } M^r_{s-1}(G \sm x) + \mbox{rk } M^{r-1}_s(G/x)\}.$$
\end{lemma}

\nib{Proof.} First we note an identity for inclusion matrices.
Suppose $H$ is a $t$-graph, $u \le t$. Let $K$ be the complete
$u$-graph on $V(H)$. Then $M^t_u(H)M^u_{u-1}(K)=(t-u+1)M^t_{u-1}(H)$.
To see this, note that if $A \in H$ and $|S|=u-1$ then the $(A,S)$
entry on the left hand side is either $0$ if $S \not\subset A$, or
otherwise the number of $u$-sets $U$ with $S \subset U \subset A$,
i.e. $t-u+1$, which agrees with the definition of the right hand side.

To write $M^r_s(G)$ in a convenient form we organise the rows
as $R = R_1 \cup R_2$ and columns as $C = C_1 \cup C_2$,
where $R_1$ corresponds to those sets of $G$ that contain $x$
and $C_1$ corresponds to all $s$-sets of $V(G)$ that contain $x$.
This gives the block form
$$ M^r_s(G) = \left( \begin{array}{cc}
M^{r-1}_{s-1}(G/x) & M^{r-1}_s(G/x) \\
0 & M^r_s(G \sm x)
\end{array} \right),$$
from which we obtain the first lower bound on the rank.
Let $M_1, M_2$ be the submatrices corresponding to the columns in
$C_1, C_2$ respectively
and $K$ be the complete $s$-graph on $V(G) \sm \{x\}$.
Now we apply the row and column operations
$$M_1' = (s-r)(r-s+1)^{-1}(M_1 - (r-s)^{-1}M_2M^s_{s-1}(K)).$$
Since $M^{r-1}_s(G/x) M^s_{s-1}(K) = (r-s) M^{r-1}_{s-1}(G/x)$
and $M^r_s(G \sm x) M^s_{s-1}(K) = (r-s+1) M^r_{s-1}(G \sm x)$
we obtain a matrix with block form
$$ \left( \begin{array}{cc}
0 & M^{r-1}_s(G/x) \\
M^r_{s-1}(G \sm x) & M^r_s(G \sm x)
\end{array} \right),$$
which gives the second lower bound on the rank. \qed

\medskip

Given this recursion, it is natural to think that for a general
size $|G|$ of the $r$-graph $G$ it may be optimal to take an initial
segment of the colex order. To explain this point further
we will briefly describe some properties of the order,
and we refer the reader to the survey \cite{FT} for
more information. Write $|G|$ in cascade form: the unique
expression
$|G| = \binom{n_r}{r} + \binom{n_{r-1}}{r-1} + \cdots +
\binom{n_j}{j}$ where $n_r > n_{r-1} > \cdots > n_j \ge j \ge 1$.
Using the natural numbers as our underlying ordered set, the
initial segment of size $G$ consists of all $r$-subsets
of $[n_r]$, all $r$-sets obtained by adding $n_r+1$ to
an $(r-1)$-subset of $[n_{r-1}]$, ... , and all $r$-sets
obtained by adding $n_r+1, n_{r-1}+1, \cdots, n_{j+1}+1$
to a $j$-subset of $[n_j]$. The shadow of this system
is the initial segment of the colex order on $(r-1)$-sets
of length $|\partial G| = \binom{n_r}{r-1} + \binom{n_{r-1}}{r-2}
+ \cdots + \binom{n_j}{j-1}$, and the Kruskal-Katona theorem
states that this is the best possible lower bound.
Iterating, we obtain that for $s \le r$ the $s$-shadow is the initial
segment of the colex order on $s$-sets of length
$|\partial^r_s G| = \binom{n_r}{s} + \binom{n_{r-1}}{s-1}
+ \cdots + \binom{n_j}{j-r+s}$, where $\binom{m}{i}$ is defined
to be zero for $i<0$.
Considering the decomposition used in Lemma \ref{rk-induct}, with
$x=n_r+1$, it is not hard to see that
$\mbox{rk } M^r_s(G) = \mbox{rk } M^r_s(G \sm x) +
\mbox{rk } M^{r-1}_{s-1}(G/x) = \binom{n_r}{s} +
 \mbox{rk } M^{r-1}_{s-1}(G/x)$ (using Gottlieb's Theorem),
so iterating we obtain $\mbox{rk } M^r_s(G) = |\partial^r_s G|$.

However, the rank of $\mbox{rk } M^r_s(G)$ may not be as large
as the $s$-shadow. For example consider the $2$-graph $C_4$
(a $4$-cycle). The size of its shadow is $4$, which is as small
as possible for a graph with $4$ edges, but its inclusion
matrix $M^2_1(C_4)$ has rank $3$. This is not merely an effect
for `small numbers' as we can use it as a building
block in larger examples: pick a number $n>5$ and
consider the $3$-graph $K^3_{n-1} \cup \{ 12n, 23n, 34n, 14n\}$.
Thus there is no direct algebraic analogue of the Kruskal-Katona
theorem. There is an an algebraic analogue of Lov\'asz's theorem,
at least for large $r$-graphs, and that is the content
of Theorem \ref{kk-alg}, which we will soon prove.

First we need the following lemma, which expresses a rigidity
property of $M^r_s(K^r_n)$ that seems independently interesting.

\medskip

\nib{Lemma \ref{full-rk}.}
Suppose $0 \le s \le r < n/2$ and $G = K^r_n \sm F$ is an $r$-graph
on $[n]$ with $|F| < \binom{r}{s}^{-1} \binom{n}{r-s}$. Then
$\mbox{rk } M^r_s(G) = \binom{n}{s}$.

\medskip

\nib{Proof.}
We argue by induction on $s$ and $r-s$. The two base cases are
straightforward: if $r=s$ then $|F|<1$, so $G=K^r_n$ and
$\mbox{rk } M^r_s(K^r_n) = \binom{n}{s}$ by Gottlieb's Theorem
(since $n > 2r$);
if $s=0$ then $|G|>0$ and $\mbox{rk } M^r_0(G)=1$.
For the induction step we choose a vertex $x$ of minimum degree
in $F$, so that
$$d_F(x) \le n^{-1} \sum_{v \in V(G)} d_F(v) =
r|F|/n < \frac{r}{n} \binom{r}{s}^{-1} \binom{n}{r-s} =
\binom{r-1}{s}^{-1} \binom{n-1}{r-s-1}.$$
By relabelling we can assuming that $x=n$.
Now $G/x = K^{r-1}_{n-1} \sm (F/x)$ and $|F/x| = d_F(x)$,
so by induction hypothesis
$\mbox{rk } M^{r-1}_s(G/x) = \binom{n-1}{s}$.
Also $G\sm x = K^r_{n-1} \sm (F\sm x)$
and $|F\sm x| \le |F| <  \binom{r}{s}^{-1} \binom{n}{r-s}
< \binom{r}{s-1}^{-1} \binom{n-1}{r-s+1}$ (since $n>2r$),
so by induction hypothesis
$\mbox{rk } M^r_{s-1}(G\sm x) = \binom{n-1}{s-1}$.
By Lemma \ref{rk-induct} we have
$\mbox{rk } M^r_s(G) \ge \mbox{rk } M^r_{s-1}(G \sm x)
+ \mbox{rk } M^{r-1}_s(G/x) = \binom{n-1}{s-1} + \binom{n-1}{s}
= \binom{n}{s}$. \qed

\medskip

Now we prove Theorem \ref{kk-alg}, which is as follows:
For every $r \ge s \ge 0$ there is a number $n_{r,s}$ so that
if $G$ is an $r$-graph with $|G| = \binom{x}{r} \ge n_{r,s}$ then
$\mbox{rk } M^r_s(G) \ge \binom{x}{s}$. Also, if
$r>s>0$ then equality holds only
if $x$ is an integer and $G=K^r_x$.

\nib{Proof.} We argue by induction on $r$ and $s$.
The cases $s=0$ and $r=s$
are trivial, so suppose $r>s>0$.
Suppose that $G \ne K^r_x$ is an $r$-graph
with $|G| = \binom{x}{r}$ and
$\mbox{rk } M^r_s(G) = \binom{x}{s}-h$ with $h \ge 0$.
We will show that if $n_{r-1,s-1} \ge \binom{u_0(r,s)}{r}$
(where $u_0(r,s)$ is given by Lemma \ref{bin-shadow}) and
$|G| \ge n_{r-1,s-1}$ then
there is some vertex $v$ so that
$G \sm v$ is an $r$-graph with $\binom{z}{r}$ edges
and $\mbox{rk } M^r_s(G \sm v) = \binom{z}{s}-h-h'$,
where $z > x-1$ and $h' > (3r)^{-r}x^{-1}$. Then we can iterate
this fact to obtain an $r$-graph $G_0$,
such that $|G_0| = \binom{z_0}{r}$ with
$n_{r-1,s-1} < |G_0| < 2n_{r-1,s-1}$
and
$$\mbox{rk } M^r_s(G_0) < \binom{z_0}{s} -
(3r)^{-r} \sum_{i=z_0+1}^x 1/i < 2n_{r-1,s-1} -
(3r)^{-r} \log \frac{x+1}{z_0+1},$$
using the estimate $\sum_{i=z_0+1}^x 1/i > \int_{z_0+1}^{x+1} dt/t =
\log \frac{x+1}{z_0+1}$.
This is less than $0$ if we suppose that
$|G| = \binom{x}{r} \ge n_{r,s}$ with $n_{r,s}$ sufficiently large, so
we will have a contradiction
to the existence of such $G$, which is the required result.

Now we show how to find the vertex $v$. We claim that
there is a vertex $v$ with
$1 \le d(v) \le \binom{x-1}{r-1} - \frac{1}{2r!} x^{r-s-1}$.
Write $n=|V(G)|$.
Then we can bound the minimum degree as
$\delta(G) \le r|G|/n = \frac{x}{n} \binom{x-1}{r-1}$.
If $n \ge x+1$ then we get
$\delta(G) \le \binom{x-1}{r-1} - \frac{1}{x+1} \binom{x-1}{r-1}
<  \binom{x-1}{r-1} - \frac{1}{2r!} x^{r-s-1}$ for large $x$,
so we can suppose that $n < x+1$. This rules out the
case when $x$ is an integer, as we are supposing $G \ne K^r_x$.
Therefore $n = \lceil x \rceil = x + \theta$ for some
$0 < \theta < 1$. Also, since $\mbox{rk } M^r_s(G) \le \binom{x}{s}
< \binom{n}{s}$ we have $|K^r_n \sm G| \ge \binom{r}{s}^{-1} \binom{n}{r-s}$
by Lemma \ref{full-rk}. This gives
$\binom{r}{s}^{-1} \binom{n}{r-s} \le \binom{n}{r} - \binom{n-\theta}{r} \
< \theta \binom{n}{r-1}$, so $\theta > \frac{s}{r} \binom{n-r+s}{s-1}^{-1}$.
Now
\begin{eqnarray*}
\binom{x-1}{r-1} - \delta(G) & \ge &  \binom{x-1}{r-1} -
 \frac{x}{n} \binom{x-1}{r-1}
= \frac{\theta}{n} \binom{x-1}{r-1}
> \frac{s}{rn} \binom{n-r+s}{s-1}^{-1} \binom{x-1}{r-1} \\
& = & (1+o(1)) \frac{s!}{r!} x^{r-s-1} >
\frac{1}{2r!} x^{r-s-1},
\end{eqnarray*}
for large $x$, as required.

Write $\binom{z}{r} = |G \sm v| = |G|-d(v)$. Since
$d(v) < \binom{x-1}{r-1}$ we have $z>x-1$.
We consider the cases $d(v) \le n_{r-1,s-1}$ and $d(v) > n_{r-1,s-1}$
separately. First suppose that $d(v) \le n_{r-1,s-1}$.
Then $n_{r-1,s-1} \ge \binom{x}{r} - \binom{z}{r}
> (x-z) \binom{z-1}{r-1} > (x-z)\binom{x-2}{r-1}$ and
$\binom{x}{s} - \binom{z}{s} < (x-z) \binom{x}{s-1}
< n_{r-1,s-1} \binom{x}{s-1} \binom{x-2}{r-1}^{-1} < 1/2$ for large $x$.
Since $d(v) \ge 1$ we have $\mbox{rk }M^{r-1}_{s-1}(G/v) \ge 1$,
so by Lemma \ref{rk-induct} we have
$\mbox{rk } M^r_s(G \sm v) \le \mbox{rk } M^r_s(G)
- \mbox{rk }M^{r-1}_{s-1}(G/v) \le \binom{x}{s}-h-1
< \binom{z}{s} - h - 1/2$.

Now suppose that $d(v) > n_{r-1,s-1}$. Write $d(v) = \binom{w}{r-1}$.
Then by the induction hypothesis we have
$\mbox{rk }M^{r-1}_{s-1}(G/v) \ge \binom{w}{s-1}$.
Now $\binom{x}{r} = \binom{z}{r} + \binom{w}{r-1}$,
where $\binom{w}{r-1} = d(v) < \frac{x-1}{r-1} - \frac{1}{2r!} x^{r-s-1}$,
so by Lemma \ref{bin-shadow} we have
$\binom{x}{s} = \binom{z}{s} + \binom{w}{s-1} - h'$, with
$h' > (3r)^{-r}x^{-1}$.
Now by Lemma \ref{rk-induct} we have
$\mbox{rk } M^r_s(G \sm v) \le \mbox{rk } M^r_s(G)
- \mbox{rk }M^{r-1}_{s-1}(G/v) \le \binom{x}{s} - h - \binom{w}{s-1}
= \binom{z}{s} - h - h'$, as required.

Either way we obtain the vertex $v$ required in the first
paragraph of the proof, so we are done. \qed

\section{Concluding remarks}

The argument for proving Theorem \ref{intstab2} via
expansion in the Cayley graph
applies generally to any extremal problem on $k$-graphs
with the property that
when one restricts to an interval there can be at most $k$ sets,
with equality exactly when they all contain some fixed point.
More generally, it gives a strategy to prove a stability
theorem for any extremal problem that can be uniformly covered by
`simpler instances' via the action of a group $G$, provided that we
have a characterisation of the maximum constructions for the
simpler instances and a set of generators for $G$ that is
`well behaved' with respect to the constructions and give
a Cayley graph with good expansion. Such a stability result could in
turn be used as part of the stability method for solving the
original problem (see \cite{KM} for an example of the stability method
and references to many other examples). We hope to return to this
idea in future work.

We have proved an algebraic analogue of Lov\'asz's theorem, but it
is natural to ask if there can be any algebraic analogue of
the Kruskal-Katona theorem, even though we have observed that
there are other factors to take into account, and so the
full description of the optimal constructions may be very
complicated.

Constructions of explicit rigid matrices can be used to
obtain lower bounds in various notions of complexity
used in Theoretical Computer Science (see \cite{KR,Lok,P}).
As far as we can see, our rigidity result (Lemma \ref{full-rk})
does not give any non-trivial result in this arena, but
perhaps some new ideas could turn it into a useful construction.
The lemma is not exactly tight, but it is tight up to a constant,
as may be seen by fixing some set $S$ of size $s$ and letting
$F$ consist of all $r$-sets that contain $S$. Then
$|F| = \binom{n-s}{r-s} = \Theta(n^{r-s})$ and $\partial^r_s G$
does not contain $S$, so $M^r_s(G)$ does not have
full rank.

\medskip

\nib{Acknowledgements.} The author thanks Dhruv Mubayi for helpful
conversations, Benny Sudakov for bringing
references \cite{Al} and \cite{FK} to his attention, and an
anonymous referee for their careful reading of the paper.

\appendix

\section{Proofs of binomial coefficient estimates}

This appendix contains the proofs of Lemmas \ref{bin-diff}
and \ref{bin-shadow}. First we recall an identity for
binomial coefficients (see Ex 1.42(i) in \cite{Lo}).

\begin{equation} \label{lo1}
\binom{x+y}{r} = \sum_{j=0}^r \binom{x+j-1}{j} \binom{y-j}{r-j}.
\end{equation}

Another exercise in \cite{Lo}, 1.43(e), states that
\begin{equation} \label{lo2}
\frac{d}{dx} \binom{x}{r} = \sum_{i=1}^r \frac{1}{i} \binom{x-i}{r-i}.
\end{equation}
Note that this is a strictly increasing function of $x$ for $x \ge r-1$.
We also need the Mean Value Theorem from Calculus, that
if $f(x)$ is a real differentiable function
and $a>b$ then $\frac{f(a)-f(b)}{a-b}=f'(c)$ for some $a \ge c \ge b$.
Furthermore, if $f'(t)$ is a strictly increasing function
we can take $a < c < b$.

\nib{Proof of Lemma \ref{bin-diff}}.
Suppose $x>y \ge r-1$. Write $f(z) = \binom{z}{r}$.
By the Mean Value Theorem we can write
$\frac{f(x)-f(y)}{x-y}=f'(c)$, for some $y < c<x$.
Then
$$\binom{x}{r} - \binom{y}{r} = (x-y)f'(c) > (x-y)f'(y)
 > (x-y) \binom{y-1}{r-1},$$
by equation (\ref{lo2}). Also
$$\binom{x}{r} - \binom{y}{r} = (x-y)f'(c) < (x-y)f'(x)
 \le (x-y) \sum_{i=1}^r \binom{x-i}{r-i}
  = (x-y) \binom{x}{r-1},$$
applying equation (\ref{lo1}) with $x$ replaced by $1$,
$y$ replaced by $x-1$ and $r$ replaced by $r-1$. \qed

\medskip

Next we will prove Lemma \ref{bin-shadow}, which can be
regarded as a defect form of Ex 13.31(a) in \cite{Lo}, and indeed
our proof involves a more careful analysis of what is going
on inside Lov\'asz's proof. First we need to give a separate
argument for the case $r=2$, which is easy.

\begin{lemma} \label{bin-shadow-2}
Suppose $C>0$ and $\binom{u}{2} = \binom{v}{2} + w$
with $1 \le w < u-1- C$. Then $u < v + 1 - C/u$.
\end{lemma}

\nib{Proof.} Since $w < u-1$ we have $v>u-1$. Also
$0 = v(v-1)+2w-u(u-1) = 2w - (u-v)(u+v-1)$,
so $1+v-u = 1 - \frac{2w}{u+v-1} > 1 - \frac{2(u-1-C)}{2u-2}
= C/(u-1) > C/u$. \qed

Now we can assume $r \ge 3$. We prove the following lemma, in
which statement (3) is the result we want, Lemma
\ref{bin-shadow}.

\begin{lemma} \label{bin-shadow-r}
Suppose $r \ge 3$, $\binom{u}{r} = \binom{v}{r} + \binom{w}{r-1}$
with $C' \le \binom{w}{r-1} < \binom{u-1}{r-1} - C$,
for some $C, C' > 0$ and $u$ is sufficiently large.
For $1 \le s \le r$ write
$X_s = \binom{v}{s} + \binom{w}{s-1} - \binom{u}{s}$.
Then

\noindent (1) if $w > u - 3r$ then $X_{r-1}>C/3u$,

\noindent (2) if $w < u - 2r$ then $X_{r-1} > \min\{1/4r!,C'\}$, and

\noindent (3) if $s \le r-1$, $C'=1$ and $C = \frac{1}{2r!} u^{r-s-1}$ then
$X_s > (3r)^{-r}u^{-1}$.
\end{lemma}

\nib{Proof.}
First consider the (possibly non-existent) case $C' \le \binom{w}{r-1} \le 1$,
when we have $r-2 < w \le r-1$. Now
$(\binom{u}{r} - \binom{v}{r}) - (\binom{u}{r-1} - \binom{v}{r-1})
= \binom{u-1}{r-1} - \binom{v-1}{r-1} > 0$ since $u>v>r-1$
so $X_{r-1} > \binom{w}{r-2} + \binom{v}{r} - \binom{u}{r}
= \binom{w}{r-2} - \binom{w}{r-1}
= \left( \frac{r-1}{w-r+2} - 1 \right) \binom{w}{r-1}
\ge (r-2) \binom{w}{r-1} \ge C'$.

Now we suppose that $\binom{w}{r-1} > 1$ and, following \cite{Lo},
introduce the change of variables $w = t+r-1$, $u'=u-t$, $v'=v-t$.
Note that $t>0$, $u' = (u-w)+r-1 > r$ and $v'>r-1$.
By identity (\ref{lo1}) we have
$\binom{u}{r} = \sum_{j=0}^r \binom{t+j-1}{j} \binom{u'-j}{r-j}$,
$\binom{v}{r} = \sum_{j=0}^r \binom{t+j-1}{j} \binom{v'-j}{r-j}$,
and $\binom{w}{r-1} = \sum_{j=0}^{r-1} \binom{t+j-1}{j}$,
so
\begin{equation} \label{x:r}
0 = X_r = \sum_{j=0}^{r-1} \binom{t+j-1}{j} \phi_{r-j},
\end{equation}
where $\phi_i = \binom{v'-r+i}{i} + 1 - \binom{u'-r+i}{i}$.
Similarly we have
\begin{equation} \label{x:r-1}
X_{r-1} = \sum_{j=1}^{r-1} \binom{t+j-1}{j-1} \phi_{r-j} =
\sum_{j=1}^{r-1} \frac{j}{t} \binom{t+j-1}{j} \phi_{r-j}.
\end{equation}
Now
$$\phi_{i-1} = \left( \frac{i}{u'-r+i} \right) \phi_i
+ \left(1 - \frac{i}{u'-r+i} \right) +
\left(\frac{i}{v'-r+i} - \frac{i}{u'-r+i} \right)
\binom{v'-r+i}{i},$$
so if $\phi_i>0$ we have $\phi_{i-1}>0$. Therefore there is
some $0 \le k \le r$ so that $\phi_i > 0$ for $1 \le i \le k$
and $\phi_i \le 0$ for $k < i \le r$.

Note that $\phi_1 = v'+1-u' = v+1-u>0$, so $k \ge 1$. We have
$\binom{v}{r} = \binom{u}{r} - \binom{w}{r-1}
> \binom{u}{r} - \binom{u-1}{r-1} + C = \binom{u-1}{r} + C$ so
$C < \binom{v}{r} - \binom{u-1}{r} < (v-u+1) \binom{v}{r-1}$
and $\phi_1 > C \binom{v}{r-1}^{-1}$. Let $k' = \max\{k,3/2\}$.
By equations (\ref{x:r-1}) and (\ref{x:r}) we have
\begin{equation} \label{phi1}
X_{r-1} > \frac{r-k'}{t} \sum_{j=0}^{r-1} \binom{t+j-1}{j} \phi_{r-j}
+ \frac{1}{2t} \binom{t+r-2}{r-1}\phi_1 =
\frac{1}{2t} \binom{t+r-2}{r-1}\phi_1.
\end{equation}

If $w > u - 3r$ then equation (\ref{phi1}) gives
$$X_{r-1} > \frac{1}{2u} \binom{u-3r-1}{r-1} \cdot C \binom{v}{r-1}^{-1}
= (1+o(1)) \frac{1}{2u} \cdot \frac{u^{r-1}}{(r-1)!} \cdot C
\cdot \frac{(r-1)!}{u^{r-1}} \sim C/2u > C/3u$$
for large $u$, which proves (1). Now suppose $w < u-2r$,
so $u' = u - w + r - 1 > 3r-1$ and $v'>u'-1>3r-2$.
If $\phi_1 \ge 1/2$ then equation (\ref{phi1}) gives
$X_{r-1} \ge \frac{1}{2t} \binom{t+r-2}{r-1} \cdot 1/2
= \frac{\prod_{i=1}^{r-2} (t+i)}{4(r-1)!} > 1/4r!$,
i.e. (2) holds, so we can suppose $\phi_1 < 1/2$.
Then $u'-v' = 1-\phi_1 > 1/2$ and
$1 - \phi_2 = \binom{u'-r+2}{2} - \binom{v'-r+2}{2}
= \frac{1}{2}(u'-v')(u'+v'-2r+3) > r$, so $\phi_2 < 1-r \le -2$.
Then $k=1$ and by equations (\ref{x:r-1}) and (\ref{x:r}) we have
\begin{eqnarray*}
X_{r-1} & > & \frac{r-1}{t} \sum_{j=0}^{r-1} \binom{t+j-1}{j} \phi_{r-j}
- \frac{1}{t} \binom{t+r-3}{r-2}\phi_2 =
 \frac{1}{t} \binom{t+r-3}{r-2}(-\phi_2) \\
& = & \frac{\prod_{i=1}^{r-3} (t+i)}{(r-2)!}(-\phi_2)
> 2/(r-2)! > 1/4r!,
\end{eqnarray*}
so (2) holds in either case.

Finally, suppose $C'=1$ and $C = \frac{1}{2r!} u^{r-s-1}$.
To prove (3) we consider the cases $w<u-2r$ and $w \ge u-2r$
separately. If $w<u-2r$ then for $1 \le i \le r-s$
we prove $X_{r-i} > 1/4r!$ by induction on $i$.
The base case $i=1$ holds by (2). For the induction step,
suppose $X_{r-i} > 1/4r!$ for some $i \ge 1$. Define
$u_i$, $w_i$ by $\binom{w_i}{r-i-1} = \binom{w}{r-i-1} - 1/4r!$
and $\binom{u_i}{r-i} = \binom{v}{r-i} + \binom{w_i}{r-i-1}$.
Then
$\binom{u}{r-i} < \binom{u_i}{r-i} = \binom{v}{r-i}  + \binom{w_i}{r-i-1}$,
$w_i < w < u-2r < u_i-2r$
and $\binom{w_i}{r-i-1} > 1 - 1/4r! > 1/4r!$,
so applying (2) with $r$ replaced by $r-i$ we have
\begin{eqnarray*}
X_{r-i-1} & = & \binom{v}{r-i-1} + \binom{w}{r-i-2} - \binom{u}{r-i-1} \\
& > & \binom{v}{r-i-1} + \binom{w_i}{r-i-2} - \binom{u_i}{r-i-1} > 1/4r!,
\end{eqnarray*}
as required. Since $1/4r! > (3r)^{-r}u^{-1}$, we have
proved (3) in the case $w<u-2r$.

On the other hand, if $w \ge u-2r$ then for $1 \le i \le r-s$
we prove $X_{r-i} > x_i = (2r!)^{-1} 3^{-i} u^{r-s-1-i}$ by induction on $i$.
The base case $i=1$ holds by (1), since $w \ge u-2r > u-3r$.
For the induction step
suppose $X_{r-i} > x_i$ for some $i \ge 1$.
Define $u_i$, $w_i$ by $\binom{w_i}{r-i-1} = \binom{w}{r-i-1} - x_i$
and $\binom{u_i}{r-i} = \binom{v}{r-i} + \binom{w_i}{r-i-1}$.
Then $u<u_i$ and
$\binom{u_i}{r-i} = \binom{v}{r-i} + \binom{w_i}{r-i-1}
< \binom{u}{r-i} + \binom{u}{r-i-1} = \binom{u+1}{r-i}$,
so $u_i < u+1$. Also $w-w_i = o(u)$, since $x_i = O(u^{r-i-2})$.
In fact
$x_i = \binom{w}{r-i-1} - \binom{w_i}{r-i-1} > (w-w_i)\binom{w_i}{r-i-2}$,
so
$$w-w_i < x_i \binom{w_i}{r-i-2}^{-1} < (2r!)^{-1} u^{r-s-1-i} \cdot
\frac{(r-i-2)!}{(w_i-r)^{r-i-2}} <
\left( \frac{u}{w_i-r} \right)^{r-i-2} u^{1-s}
< 1+o(1),$$
which implies $w_i>w-2$.
Therefore $w_i>w-2>u-2r-2>u_i-2r-3 \ge u_i-3r$.
Since $\binom{w_i}{r-i-1} = \binom{w}{r-i-1} - x_i
< \binom{u_i-1}{r-i-1} - x_i$ we can apply (1), or Lemma
\ref{bin-shadow-2} in the case $s=1$ and $i=r-2$,
with $r$ replaced by $r-i$ to get
\begin{eqnarray*}
X_{r-i-1} & = & \binom{v}{r-i-1} + \binom{w}{r-i-2} - \binom{u}{r-i-1} \\
&>&  \binom{v}{r-i-1} + \binom{w_i}{r-i-2} - \binom{u_i}{r-i-1}
> x_i/3u = x_{r-i-1}.
\end{eqnarray*}
Then $X_s > 3^{s-r}(2r!)^{-1}u^{-1} > (3r)^{-r}u^{-1}$, so
we have (3) in both cases, and the
lemma is proved. \qed

\end{document}